\documentclass[12pt]{article}
\pdfoutput=1
\usepackage{a4wide}
\usepackage{amsmath}
\usepackage{amssymb}
\usepackage{comment}
\usepackage{hyphenat}
\usepackage{multicol}

\usepackage{graphicx}
\usepackage{epsfig}
\usepackage{epstopdf}
\usepackage{color}
\usepackage{float}
\newtheorem{theorem}{Theorem}[section]
\newtheorem{definition}[theorem]{Definition}

\newtheorem{corollary}[theorem]{Corollary}

\newcommand{\qed}{\hfill$\Box$}

\title{Analysis of some Epidemic Models in complex networks and some ideas about isolation strategies}
\author{Carlos Rodr\'iguez Lucatero \\
Departamento de Tecnolog\'{\i}as de la Informaci\'on \\ 
Universidad Aut\'onoma Metropolitana-Cuajimalpa\\
Torre III \\
Av. Vasco de Quiroga 4871 \\
Col.Santa Fe Cuajimalpa, M\'exico, D. F.\\
C.P. 05348, M\'exico \\
email:crodriguez@correo.cua.uam.mx \\
}
\date{}
\begin{document}
\maketitle
\begin{abstract} \label{resumen}
Many models of virus propagation in Computer Networks inspired by {\bf SIS,SIR,}\\
{\bf SEIR}, etc. epidemic disease propagation mathematical models that can be found in the epidemiology field have been proposed in the last two decades. The purpose of these models has been to determine the conditions under which a virus becomes rapidly extinct in a network. 
The most common models proposed in the field of virus propagation in networks are inspired by SIS-type models or their variants. In such models, the conditions that lead to a rapid extinction of the spread of a computer virus have been calculated and its dependence on some parameters inherent to the mathematical model has been observed. In this article we will try to analyze a particular model proposed in the past and show through simulations the influence that topology has on the dynamics of the spread of a virus in different networks.
A consequence of knowing the impact of the topology of a network can serve to propose effective isolation strategies to reduce the spread of a virus through modifications to the original network of contacts. I will talk about this subject at the final section of the present article.

{\em Mathematics Subjects Classification}: 
\\
\end{abstract}

%\linenumbers

\section{Introduction} \label{Introduccion} 
The issue of virus spread as well as the conditions under which it is extinguished or contaminated most of the nodes of a computer network has been studied for at least two decades. The models that have been proposed are based on mathematical models of virus spread in the field of epidemiology \cite{Chakrabarti},\cite{Chakrabarti2},\cite{Chakrabarti3}. These models are differential equations that try to capture the dynamics of a virus spread in order to answer questions such as, How long will the epidemic disappear? At what point will it peak? Will it remain endemically at a certain level of infected? Will it generate permanent immunity to those who suffered from the disease?.\\
These mathematical models emerged from the area of epidemiology are known as compartmental because they model the process as states in which an individual finds himself during an epidemic and transition probabilities between these states. In other words, an individual can be susceptible, exposed, infected, or recovered from a virus, and these states are modelled as compartments. From these compartments and the transition probabilities between them, relationships can be established in the form of differential equations that describe the dynamics of the spread of a disease \cite{Hethcote}.
The compartments used by one of those models depend on the disease to be modelled, since there are diseases that produce permanent immunity such as measles, while some diseases such as seasonal influenza do not produce permanent immunity.
The spread of a disease in epidemic processes shares similarities with the massive attacks on a computer network and therefore some ideas from the mathematical modelling of epidemics served as the basis for models of virus spread in computer networks. However, it is worth mentioning that some problems that arise in the study of the spread of viruses in computer networks, as well as their respective solutions, may be useful for epidemiology.
There have been very famous massive attacks on computer networks. One of the first denial-of-service-type attacks turned 20 on February 7, 2020, and was conducted by a 15 years old Canadian hacker  whose pseudonym was {\em Mafiaboy} and whose real name is Michael Calce.
Denial of service type attacks consist of sending a huge number of service request packets to a target server in such a way that it exceeds the capacity of said server to respond to so many orders, thus causing it to crash. This massive attack revealed the vulnerability of networks such as the internet and led to a study of the causes of said vulnerability and and gave birth to a new type of computer virus. This vulnerability aroused the interest of researchers in the field of network security and they realized that certain topologies favoured a faster dissemination of information than others or that a certain type of interconnection kept the nodes of a network connected despite the fact that some lines were faulty. Some studies showed that the type of interconnection structures that are formed in networks such as the Internet, Facebook, tweeter have similar characteristics and that they produce the appearance of giant components as well as the phenomenon of small worlds.
The formalization of the relationship of the topology of the networks and the appearance of the aforementioned phenomena uses the theory of random graphs as a mathematical tool.
The type of interconnection structures that are formed in social networks such as facebook, tweeter and the internet have been characterized as particular graphs whose distribution of degrees of the nodes is known as the type of power laws due to the mathematical form that expresses said distribution which in algebraic terms would be $p_{k} \thicksim C k^{-\gamma} $ 
with $2 < \gamma \leq 3$. The formation of networks with these topological characteristics are closely related to the type of connection protocol used and is known as {\em preferential attachment}.
It is said that this type of network is both robust and vulnerable since if they were removed by randomly choosing more than $90\%$ of the edges of the associated graph, it would continue to be a connected graph, while if they were removed by strategically choosing a number small number
($2.3\%$) of vertices of the associated graph, the original graph would be split into unconnected components \cite{Durrett2}.
Networks that have these characteristics of node degree distribution allow phenomena such as small worlds to appear due to the fact that their diameter is generally smaller than that of other graphs with other topologies. This also has an impact on the speed with which messages are spread and that I will illustrate through simulations in the following sections of this article.

\section{Epidemiology Mathematical models} \label{modelos}
Epidemics have been with humans for a long time. One of the oldest diseases is leprosy. Another disease that devastated Europe in the Middle Ages was the Black Death. Little was known about these diseases and how to cure them. In these circumstances, it was necessary to try to understand the way in which the epidemic process was developing in order to at least stay safe from them. In the absence of objective knowledge based on science, there was a tendency to believe that such calamities were divine punishments. Perhaps since those times the most helpful strategy consisted of isolating oneself.
Some diseases such as smallpox were treated since ancient times using the variolation method, which consisted in inoculating the scales of a sick patient in a healthy individual and observing that this allowed said individual to acquire a certain type of immune protection. The knowledge of this empirical method of immunization was imported to Europe from the colonies and later served as a starting point for the development of vaccines for this disease \cite{Hethcote}.
In search of answers to questions that disturb humanity during the years, one of the first mathematical models of the spread of smallpox appeared in the year 1760 and was proposed by Daniel Bernoulli
\cite{Bernoulli1}.
Later, deterministic mathematical models of virus spread began to be developed at the beginning of the 20th century. In 1906 a discrete time model was formulated and analyzed for the measles epidemics 
\cite{Hamer1}.
In 1911 some differential equations based models for the malaria disease were formulated in \cite{Ross1}, \cite{Ross2}, \cite{Ross3} and \cite{Ross4}. 

\begin{figure}[H] 
\centering
\includegraphics[width=9.5cm, height=9.5cm]{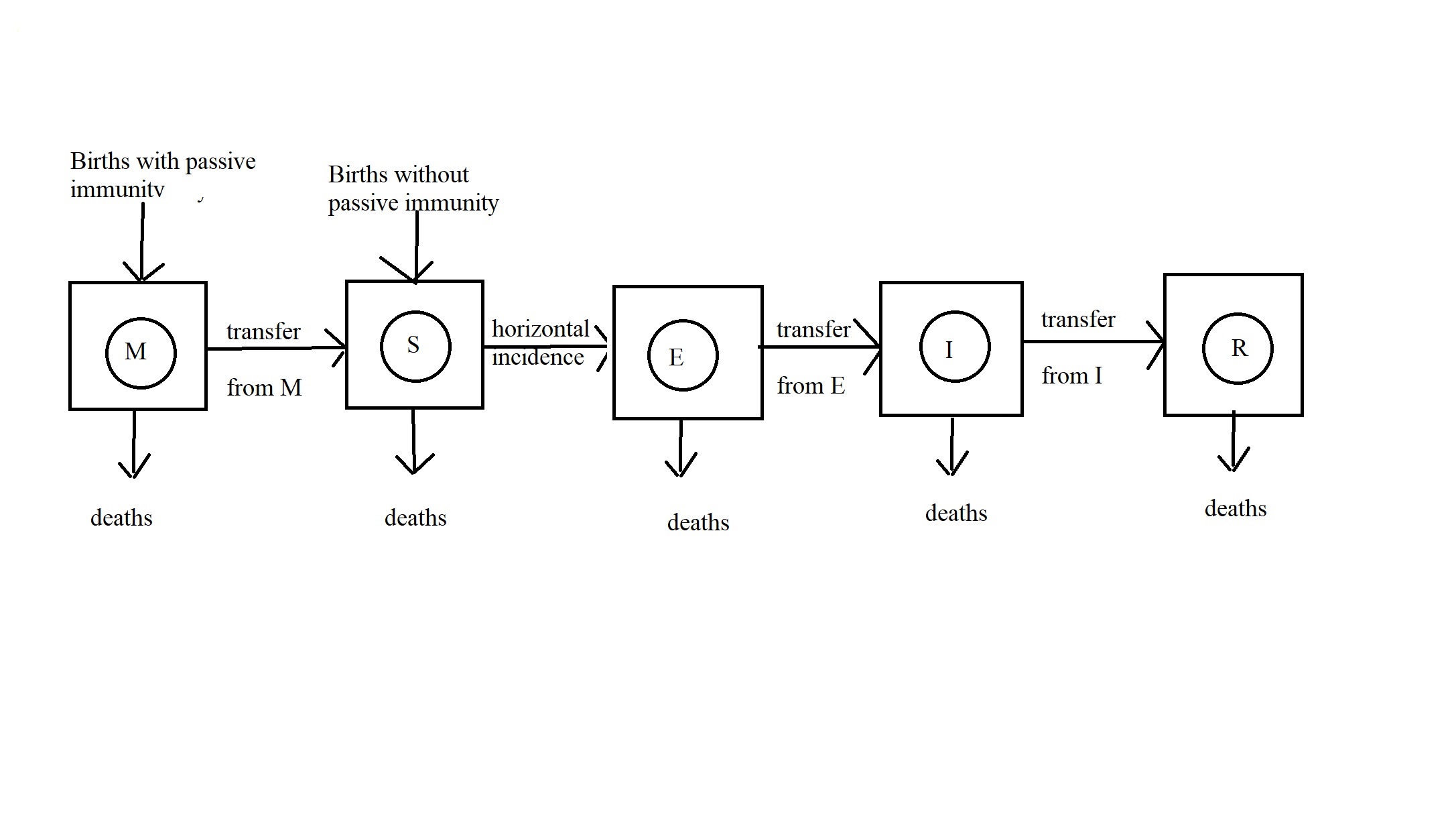}
%\centerline{\psfig{file=chakra101b.eps,width=4.2cm}}
%\centerline{\psfig{file=..\FigurasArticulo\MSEIR.jpeg,width=4.2cm}}
\vspace*{5pt}
\caption{The general transfer diagram for the MSEIR model}
\label{MSEIR}
\end{figure}

Mathematical models of this type of phenomenon clarify which are the important parameters to take into account to obtain concepts such as thresholds, basic reproduction numbers, number of contacts and replacement numbers. This in turn allows us to make computer  simulations.
Having the information provided by epidemiological models helps to know what data must be collected in an epidemic, identify trends, make estimates and calculate uncertainties in these estimates.
Later, in 1926, models were formulated from which thresholds could be calculated from which an epidemic outbreak appears and this happens when the number of susceptible reaches a critical value  \cite{Ross2},\cite{McKendrick1} and \cite{McKendrick2}. 
The model presented in \cite{McKendrick2} is a reference in the modelling with non-linear dynamic systems of the phenomenon of virus propagation in networks.
This type of model is known as composed by compartments since it conceives the states through which an individual passes during an epidemic as boxes or compartments. The compartments are labelled by the letters $M, S, E, I, R $ and $S$ that cover the different characteristics that one may have, such as having generated antibodies from the mother's womb by maternal transmission, being susceptible, being exposed, being infected and having recovered from the disease respectively.
It is entered vertically either to the state $M$ or to the state $S$ and it is exited vertically when passing to a state of death. There are horizontal transfers between state $M$ to $S$, from $S$ to 
$E$, from $E$ to $I$ and from $I$ to $R$ as it is shown in figure \ref{MSEIR}.
Some of these compartments will be present in a specific model and others will not be, depending on the particular characteristics of the disease to be modelled.
For instance, if the disease to be modelled produces permanent immunity, do not have a passive immunity transmitted by the mother and there is no incubation period, then the compartments that will be present are $S,I$ $R$ (Susceptible, Infected and Recovered) and will be called a $SIR$-model.

%\begin{figure}[th]

Frequently in epidemic models the average number of secondary infections from one individual is called the reproduction number and is denoted as $R_{0}$ and normally is bigger than $1$. When the time span of the process is quite long, the model is called {\em endemic} otherwise it is called {\em epidemic}. 
In the case of {\em endemic} models, factors as population growth or population decrease must be included and it affects the calculation of $R_{0}$. In other models the the demographic age structure are taken into account and this can change the calculation of $R_{0}$ as well.
In the next subsection I will formulate the classic mathematical epidemic models  known as 
{\bf SIR} and {\bf SIS} and show how they work with some simple simulations that I have implemented in {\bf MATLAB}.

%% Hasta aqui voy CRL 27/dic/2020

\subsection{Formulation of two well known epidemiology models} \label{SIRSIS}

The first epidemic model that I will describe is the classic {\bf SIR} model that models diseases where some individuals start out being susceptible. These individuals can transition to a state of infection by contagion and after a time they transition to a state of recovery obtaining permanent immunity.
The arrow labelled as {\em horizontal incidence} in figure \ref{MSEIR} represent the infection rate with which a susceptible individual is in contact with infected nodes and make them transit to an infected state.
$S(t)$ represent the number of susceptible in time $t$, $I(t)$ the number of
infected in time $t$, $N$ the total size of the population, $s(t)=\frac{S(t)}{N}$ the fraction of susceptible of the total population  in time $t$, $i(t)=\frac{I(t)}{N}$ the fraction of infected of the total population in time $t$, $N$ the population size,$s(t)=\frac{S(t)}{N}$ the fraction of susceptible, $i(t)=\frac{I(t)}{N}$ the fraction of infected of the total population in time $t$ and $\beta$ the average number of adequate contact or sufficient for transmission.
The average number of contacts with infected per unit time of one susceptible 
is expressed as $\frac{\beta I}{N}=\beta i$ and $\frac{\beta I}{N} S =\beta N i s$ the number of new cases per unit time. The transitions from the boxes $M,E$ and $I$ in figure \ref{MSEIR} are calculated $\delta M, \epsilon E$ and $\gamma I$. These terms represent the exponentially distributed waiting time in each box that in the cases of the $I$ compartment in figure \ref{MSEIR} the transfer rate 
$\gamma I$  $P(t)= e^{-\gamma t}$ is the fraction that is still
in the infective class $t$ units after entering this class and $\frac{1}{\gamma}$ is the mean waiting time.

The parameters $R_{0}, \sigma$ and the replacement number $R$ are related with the threshold.
For more details about the classical {\bf SIR} model consult the article \cite{Hethcote}

After having defined the relevant elements of the phenomenon of the spread of a disease, we are able to define the following $SIR$ epidemic model:

\begin{equation} \label{SIRepidemico1}
\begin{array}{ll}
      \frac{dS}{dt}=-\beta \frac{IS}{N},           & S(0)=S_{0} \geq 0\\
      \frac{dI}{dt}=\beta \frac{IS}{N} -\gamma I,  & I(0)=I_{0} \geq 0\\
      \frac{dR}{dt}=\gamma I,                      & R(0)=R_{0} \geq 0
\end{array}      
\end{equation}

where $S(t)+I(t)+R(t)=N$. If we divide by $N$ (the total population) the equations \ref{SIRepidemico1} we get

\begin{equation} \label{SIRepidemico2}
\begin{array}{ll}
      \frac{ds}{dt}=-\beta is,           & s(0)=s_{0} \geq 0\\
      \frac{di}{dt}=\beta is -\gamma i,  & i(0)=i_{0} \geq 0
\end{array}      
\end{equation}
 
with $r(t)=1-s(t)-i(t)$ where $s(t), i(t)$ and $r(t)$ are de fractions in the classes.

We are also able to define the following $SIR$ endemic model

\begin{equation} \label{SIRendemico1}
\begin{array}{ll}
      \frac{dS}{dt}= \mu N - \mu S -\beta \frac{IS}{N},           & S(0)=S_{0} \geq 0\\
      \frac{dI}{dt}= \beta \frac{IS}{N} -\gamma I - \mu I,  & I(0)=I_{0} \geq 0\\
      \frac{dR}{dt}=\gamma I - \mu R,                      & R(0)=R_{0} \geq 0
\end{array}      
\end{equation}

where $S(t)+I(t)+R(t)=N$. The $SIR$ model \ref{SIRendemico1} is almost the same as the
epidemic version \ref{SIRepidemico1} except that it has an inflow of newborns into the susceptible class at rate $\mu N$ and deaths in the classes at rates $\mu S, \mu I$ 
and $\mu R$. 
If we divide by $N$ (the total population) the equations \ref{SIRendemico1} we get

\begin{equation} \label{SIRendemico2}
\begin{array}{ll}
      \frac{ds}{dt}=-\beta is + \mu - \mu s,           & s(0)=s_{0} \geq 0\\
      \frac{di}{dt}=\beta is -(\gamma + \mu) i,  & i(0)=i_{0} \geq 0
\end{array}      
\end{equation}
 
with $r(t)=1-s(t)-i(t)$ where $s(t), i(t)$ and $r(t)$ are de fractions in the classes.

For a deep exposition of more sophisticated models and analysis of their respective thresholds consult \cite{Hethcote}. 

The next figure correspond to a simulation of the {\bf SIR} model 
%\begin{figure}[th]
\begin{figure}[H]
\centering
\includegraphics[width=15cm, height=8cm]{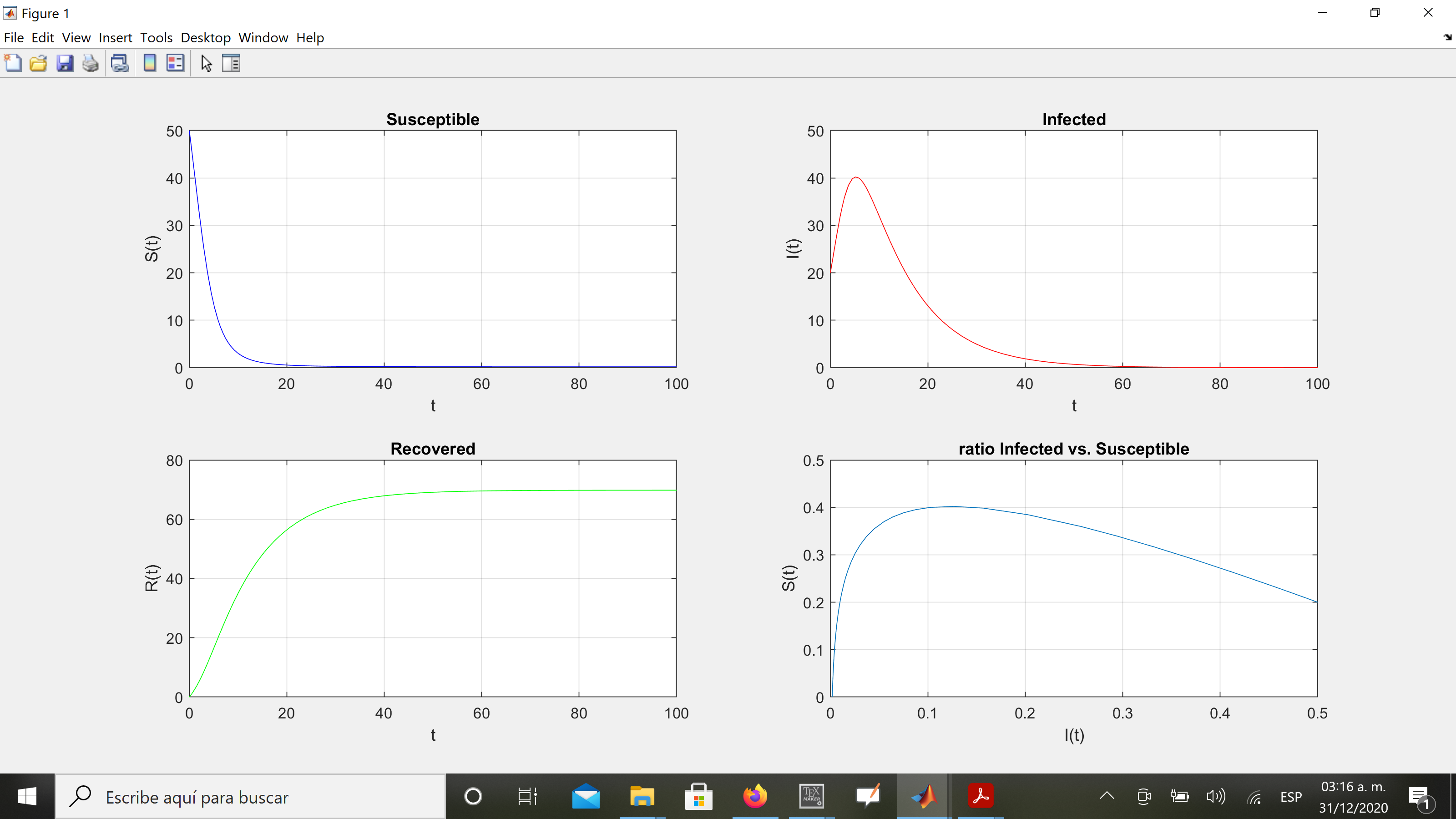}
%\centerline{\psfig{file=chakra101b.eps,width=4.2cm}}
%\centerline{\psfig{file=..\FigurasArticulo\MSEIR.jpeg,width=4.2cm}}
\vspace*{3pt}
\caption{Simulation of classical epidemic SIR model with $\beta=0.8, \gamma=0.1$}
\label{flujofaseSIR}
\end{figure}
 
 In reference to figure \ref{flujofaseSIR}, as can be seen in the graphs generated by the simulation of the classical epidemic {\bf SIR} model, the behavior in time of the infected curve reaches its peak or acme and then decreases until it is extinguished as expected, since the type of epidemic processes that the model describes is related to infections that produce permanent immunity, which is also observed in the corresponding growth in behavior over time of the recovered curve. It is also important to note that the curve of the proportion of Infected vs. Susceptibles converges to a fixed point in the coordinates $ (0,0) $, that is, the origin of said coordinate system.
 
 The classical {\bf SIS} model can be defined  as follows
 
 \begin{equation} \label{SISepidemico1}
\begin{array}{ll}
      \frac{dS}{dt}=-\beta \frac{IS}{N}+ \gamma I,           & S(0)=S_{0} \geq 0\\
      \frac{dI}{dt}=\beta \frac{IS}{N} -\gamma I,  & I(0)=I_{0} \geq 0
\end{array}      
\end{equation}
 
The next figure correspond to a simulation of the {\bf SIS} model 
%\begin{figure}[th]
\begin{figure}[H]
\centering
\includegraphics[width=15cm, height=8cm]{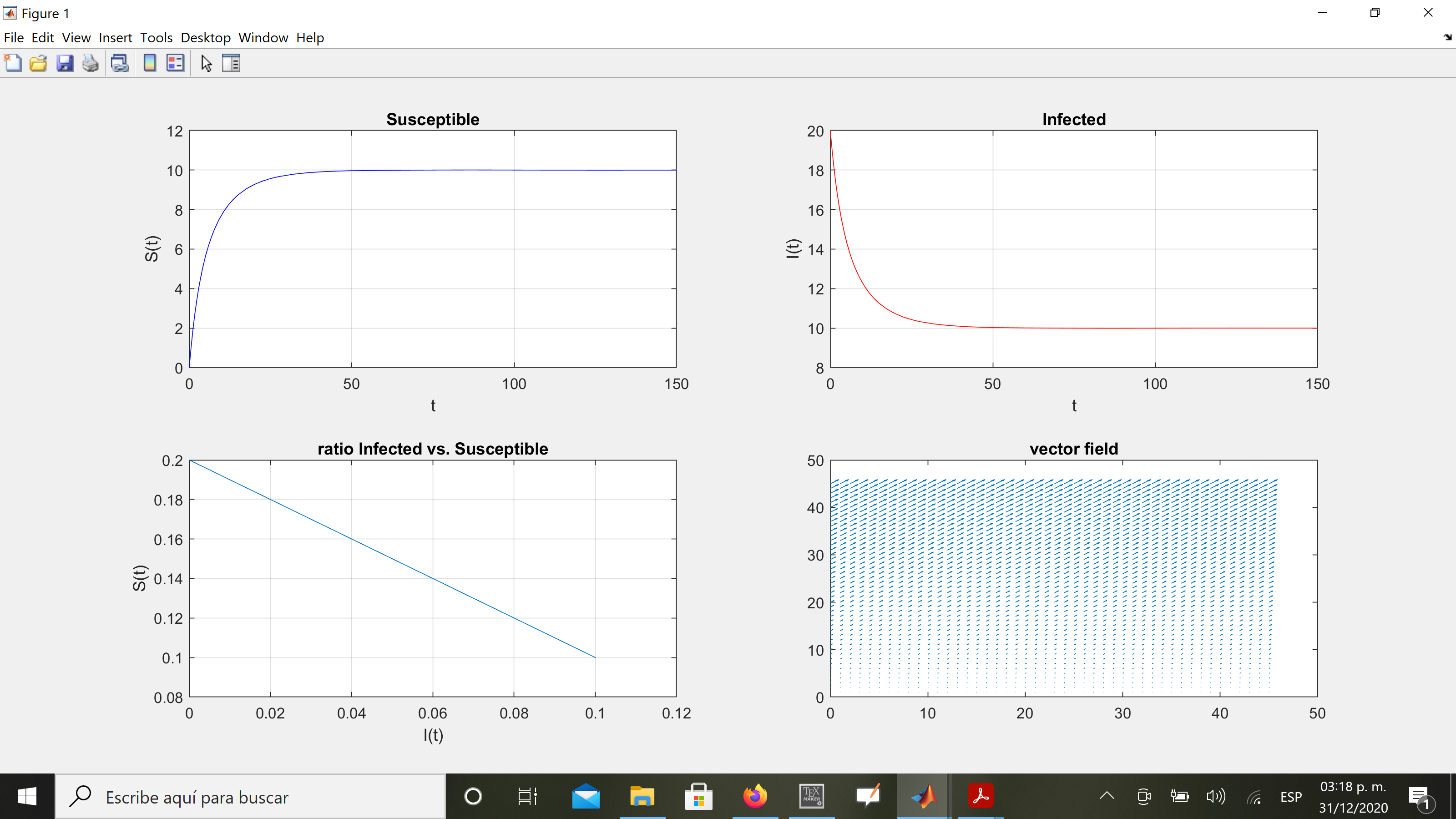}
%\centerline{\psfig{file=chakra101b.eps,width=4.2cm}}
%\centerline{\psfig{file=..\FigurasArticulo\MSEIR.jpeg,width=4.2cm}}
\vspace*{3pt}
\caption{Simulation of classical epidemic SIS model with $\beta=1, \gamma=0.1$}
\label{flujofaseSIS}
\end{figure} 

In reference to figure \ref{flujofaseSIS}, as can be seen in the graphs generated by the simulation of the classical epidemic {\bf SIS} model, the behavior in time of the infected decreases until it converges to some given level and simultaneously the number of susceptible grow until some level and stays there as expected, since the type of epidemic processes that the model describes is related to infections that produce  temporary immunity. It is also important to note that the curve of the proportion of Infected vs. Susceptible converges to a fixed point in the coordinates $ (0.1,0.1) $.

\section{Discrete epidemic models in Networks} \label{Chakrabarti}
Mathematical epidemic models assume that each individual has on average the same number of contacts. In my opinion, this hypothesis is not always fulfilled because it does not correspond to the type of structures that arise in networks such as the internet in media such as Facebook or tweeter. In addition, the type of structures that appear in these social networks partly reflects the way in which people build networks of collaboration and interaction in real life. For this reason, I believe that just as the mathematical models of epidemics have contributed to the development of discrete mathematical models to study the phenomena of virus propagation in computer networks, also the mathematical models of epidemics can benefit from the advances that take place in the study of virus spread in computer networks.
Since computer networks have grown rapidly, security problems have also increased in the same proportion. In the same way, the number of types of services offered on networks such as the Internet has grown. One service that appeared has been the distribution of content over the internet, which requires ensuring that the information reaches its destination quickly and with a good level of quality. In this type of services, P2P-type networks have been studied.
Another type of service that arises with the appearance of networks is that of sensor networks. When trying to solve both security and quality of service problems, it is necessary to resort to mathematical tools that allow modeling the phenomena inherent in computer networks and answering questions such as which network structure is the most appropriate to ensure that the network stay connected in the presence of faults in the lines or what type of topology ensures that the distance traveled by an information packet is short enough for it to reach its destination without delays.
On the other hand, the security of a network must be guaranteed and therefore it is important to know which network topology facilitates or inhibits the spread of a virus in it.
These topics have aroused the interest of researchers from fields as varied as statistical physics or experts in the field of random graphs.

The type of graphs that appear most frequently in the study of virus propagation in complex networks and that we will use to illustrate the operation of discrete SIS models are those of {\em binomial distribution of degrees, distribution of degrees in power laws, exponential type} and {\em lattices}.
For the sake of clarity I give below the definitions of these types of graphs.
\begin{definition} \label{Powerlawdef}
{\bf Power law or scale-free degree distribution graph}
Is a graph whose degree distribution of nodes follows asymptotically a power law.
More formally let $P(k)$ the fraction of the total number of nodes in a given graph that have $k$ connections with other nodes. This fraction of nodes have the following
behaviour
\begin{equation}
 P(k) \backsim k^{-\gamma}
\end{equation}

where $\gamma$ is a parameter in the interval $2 < \gamma < 3$
\end{definition}

\begin{definition}\label{Binomialdef}
{\bf Binomial degree distribution graph (Erdős–Rényi model, Barabasi-Albert model, etc.)}
Is a random graph whose degree distribution of nodes follows a binomial probability distribution law of degrees $k$ that can be formally defined as follows.
Each of the $n$ nodes of the graph is independently connected with other node with probability $p$ or not connected with probability $(1-p)$.
Let $P(k)$ the fraction of the total number of nodes in a given graph that have $k$ connections with other nodes. This fraction of nodes have the following behaviour

\begin{equation}
 P(k) = \binom{n-1}{k} p^{k}(1-p)^{n-1-k}
\end{equation}

\end{definition}

\begin{definition} \label{Exponentialdef}
{\bf Exponential degree distribution graph}
Is a random graph whose degree distribution of nodes follows a binomial probability distribution law of degrees $k$ that can be formally defined as follows.
Let $P(k)$ the fraction of the total number of nodes in a given graph that have $k$ connections with other nodes. This fraction of nodes have the following behaviour

\begin{equation}
 P(k,\lambda) = \left\{
 \begin{array}{ll}
    \lambda e^{-\lambda k} & k \geq 0 \\
    0                      & k < 0
    \end{array} \right.
\end{equation}

where $\lambda > 0$ is a parameter of the distribution called {\em rate parameter}. 

\end{definition}

\begin{definition}\label{Lattice4def}
{\bf Lattice 4 connected graph (grid graph, mesh graph, etc)}
Is a graph that each node is connected to four other nodes for all the $n$ nodes belonging to the graph.
\end{definition}

Next I will show some examples of these graphs.
\begin{figure}[H] 
\centering
\includegraphics[width=7cm, height=5cm]{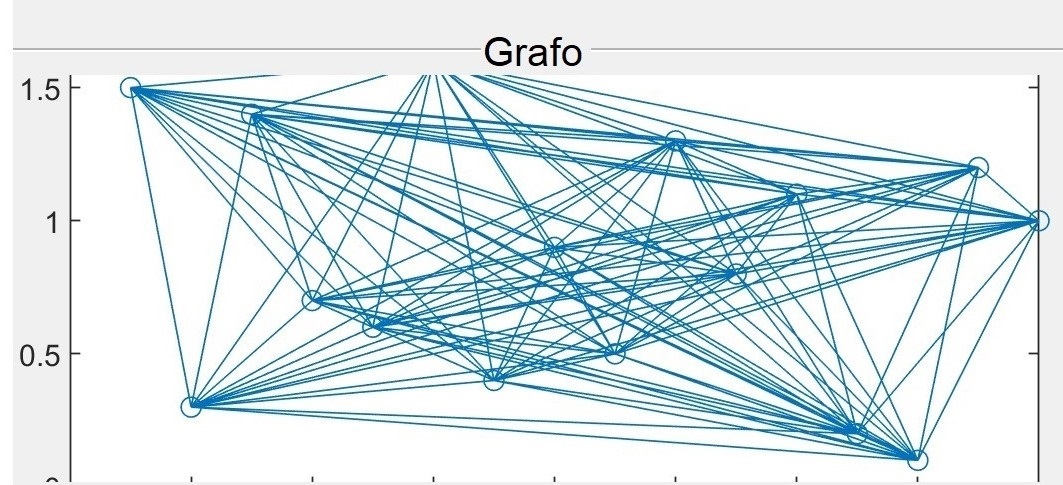}
%\centerline{\psfig{file=chakra101b.eps,width=4.2cm}}
%\centerline{\psfig{file=..\FigurasArticulo\MSEIR.jpeg,width=4.2cm}}
\vspace*{5pt}
\caption{Binomial degree distribution graph}
\label{grafobinom}
\end{figure}

\begin{figure}[H] 
\centering
\includegraphics[width=7cm, height=5cm]{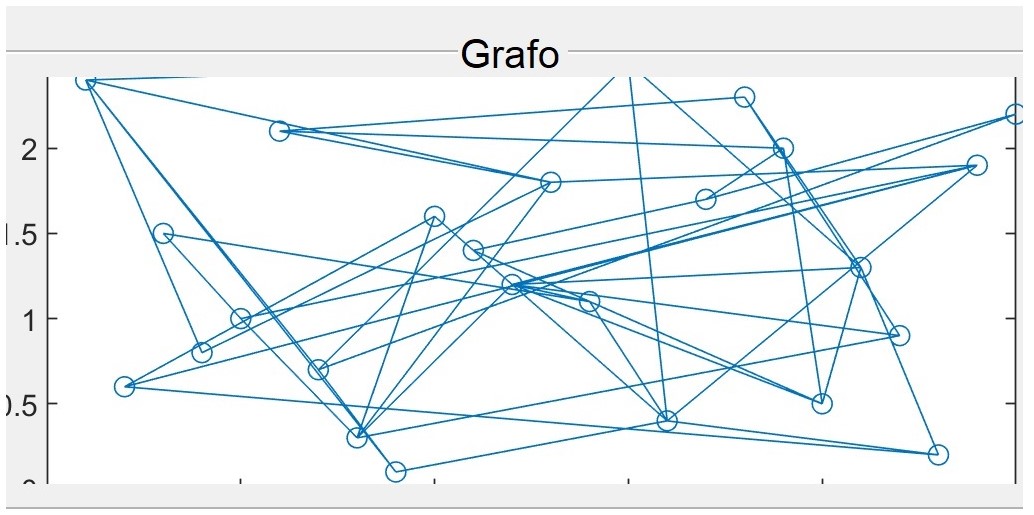}
%\centerline{\psfig{file=chakra101b.eps,width=4.2cm}}
%\centerline{\psfig{file=..\FigurasArticulo\MSEIR.jpeg,width=4.2cm}}
\vspace*{5pt}
\caption{Exponencial degree distribution graph}
\label{grafoexpon}
\end{figure}

\begin{figure}[H] 
\centering
\includegraphics[width=7cm, height=5cm]{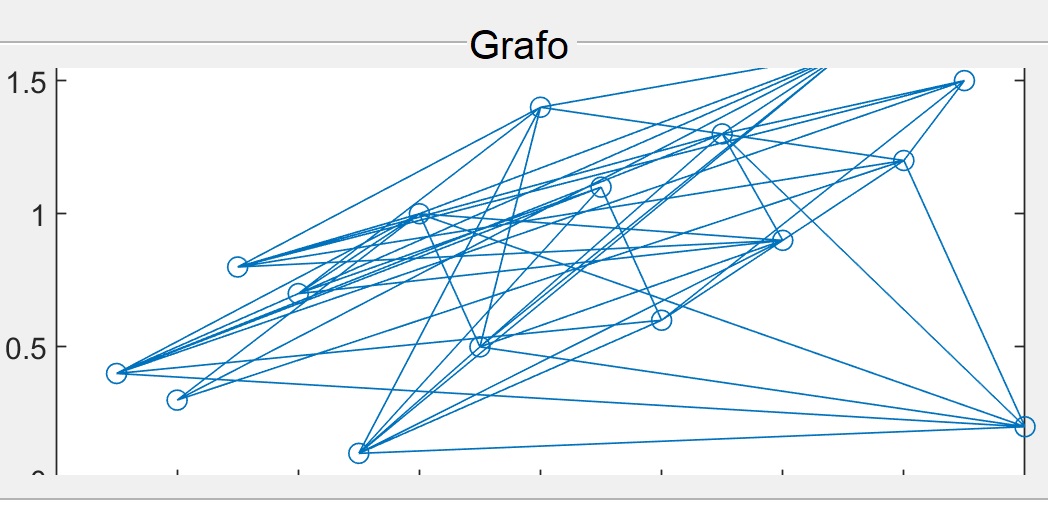}
%\centerline{\psfig{file=chakra101b.eps,width=4.2cm}}
%\centerline{\psfig{file=..\FigurasArticulo\MSEIR.jpeg,width=4.2cm}}
\vspace*{5pt}
\caption{Power law degree distribution graph}
\label{grafopower}
\end{figure}

\begin{figure}[H] 
\centering
\includegraphics[width=7cm, height=5cm]{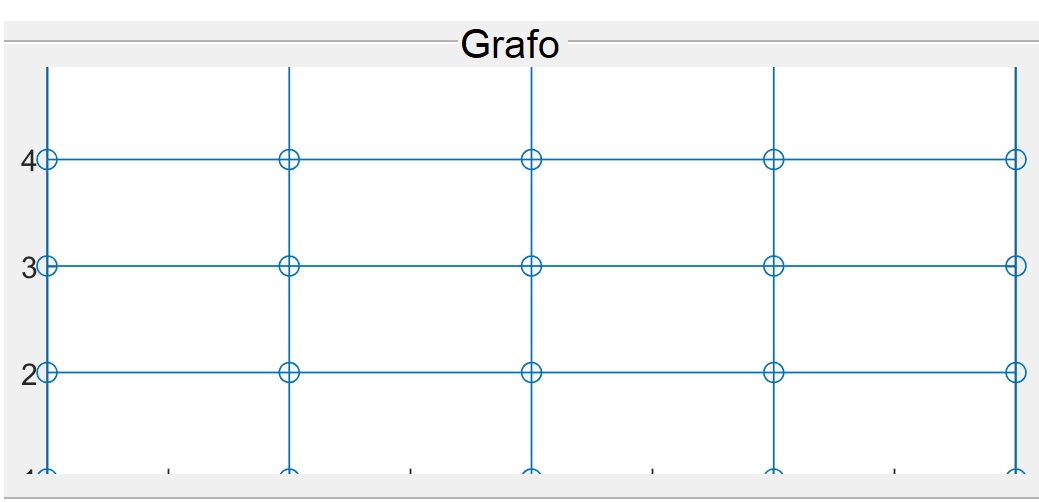}
%\centerline{\psfig{file=chakra101b.eps,width=4.2cm}}
%\centerline{\psfig{file=..\FigurasArticulo\MSEIR.jpeg,width=4.2cm}}
\vspace*{5pt}
\caption{Lattice 4 graph}
\label{grafopower}
\end{figure}

The understanding of the emergence of {\em giant components} or phenomena of {\em small worlds} that occurs on the Internet requires mathematical tools such as the theory of contact processes as well as the theory of random graphs to be able to analyze these phenomena.
 
Many research articles have been written about such subjects \cite{Durrett1},\cite{Barab0,Barab1,Barab2}, \cite{Satorras1,Satorras2,Satorras3}. Some papers about rumours spreading
on networks under the approach of contact processes have also  been written \cite{Kempe1}.
The problem of vaccine distribution on networks can be consulted in \cite{chayes1}.

Finally the use of tools such as non-linear dynamical systems, fix-point theorems for obtaining fast extinction conditions of a virus in a network combined with a discrete version of a {\em SIS} epidemic model can be found in \cite{Deepa1}. This is the model that I will describe in the following subsection.

\subsection{Discrete SIS epidemic model}
\begin{figure}[H] 
\centering
\includegraphics[width=10cm, height=7cm]{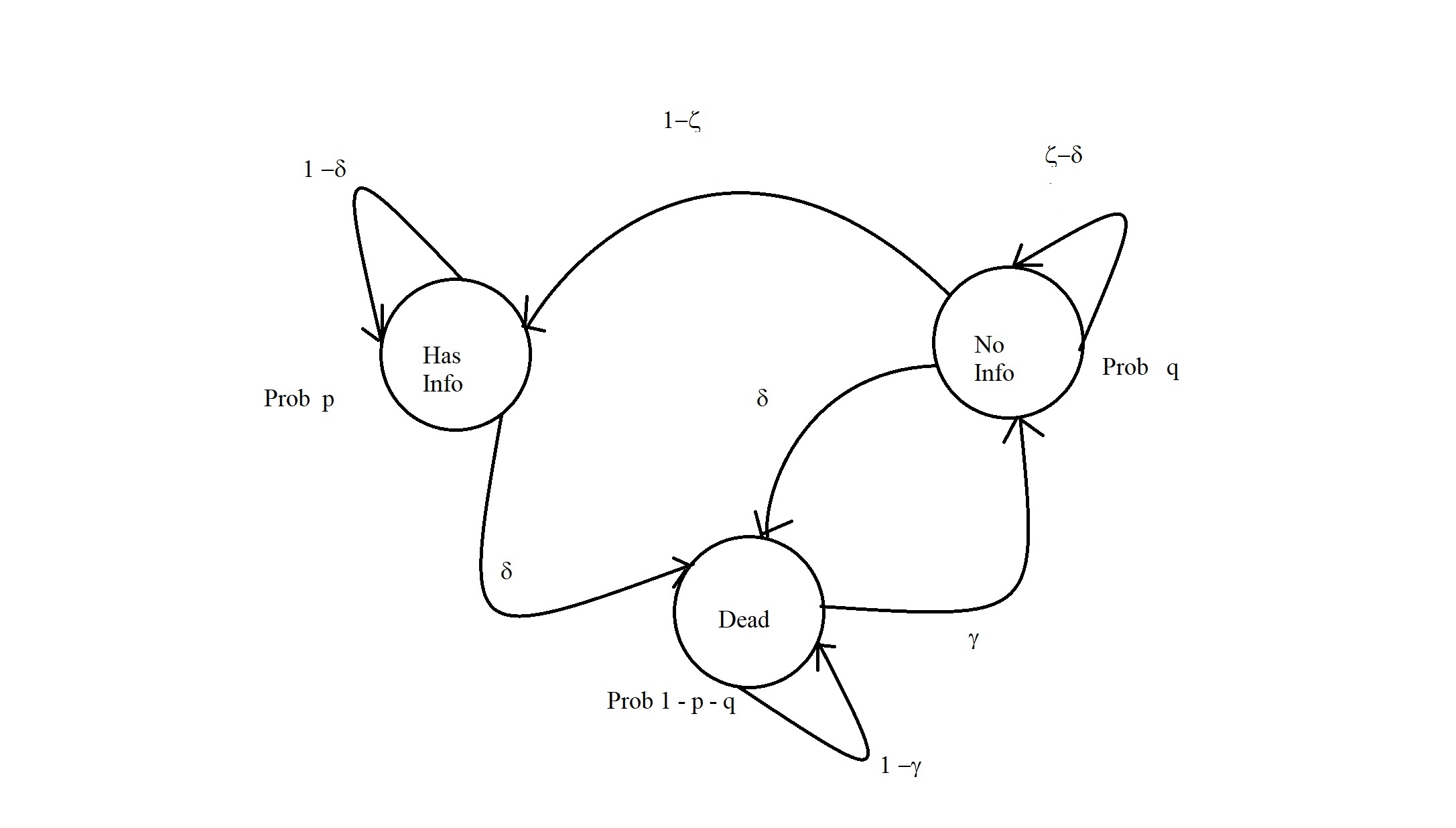}
%\centerline{\psfig{file=chakra101b.eps,width=4.2cm}}
%\centerline{\psfig{file=..\FigurasArticulo\MSEIR.jpeg,width=4.2cm}}
\vspace*{5pt}
\caption{Chakrabarti SIS model}
\label{SISChdiag}
\end{figure}

In order to understand the spread of viruses on a network,
the model proposed in \cite{Deepa1} assumes that the nodes behave according to a SIS-type model and that they are interconnected by a network.
They also assume that we take very small discrete timesteps of size $\Delta t$ where $\Delta t \rightarrow 0$. The survivability results in \cite{Deepa1} apply equally well to continuous systems. Within a $\Delta t$ time interval, each node $i$ has probability $r_{i}$ of trying to broadcast its information every time step, and each link $i \rightarrow j$ has a probability $\beta_{i,j}$ of being \emph{up}, and thus correctly propagating the information to node $j$. Each node $i$ also has a node failure probability $\delta_{i} > 0$. Every dead node $j$ has a rate $\gamma_{j}$ of returning to the \emph{up} state, but without any information in its memory. The details about the parameters of this model as well as the fast extinction conditions and stability results can be consulted in \cite{Deepa1}.
the authors of \cite{Deepa1} chose to use the non-linear dynamic systems approach and fixed point theorems. The state transitions at each node are shown in the figure \ref{SISChdiag}.
 
%\begin{figure}[th]

The node state {\em has info} can be think as being infected.
The authors of \cite{Deepa1} proposed to obtain an approximation of the threshold by describing the problem as non-linear dynamic system with $N$ variables representing the nodes and assumed that the state of two different nodes are independent.
The independence condition can be formally expressed ad follows:

\begin{equation} \label{indep}
	\zeta_{i}(t)=\prod_{j=1}^{N} (1-r_{j}\beta_{ji}p_{j}(t-1))
\end{equation}
Then equations describing the state transitions 
in the dynamic systems for each node, taking into account what is depicted in
the figure \ref{SISChdiag}, can be expressed as

\begin{equation}\label{probp}
	p_{i}(t)=p_{i}(t-1)(1-\delta_{i})+q_{i}(t-1)(1-\zeta_{i}(t))
\end{equation}

\begin{equation}\label{probq}
	q_{i}(t)=q_{i}(t-1)(\zeta_{i}(t)-\delta_{i})+(1-p_{i}(t-1)-q_{i}(t-1))\gamma_{i}
\end{equation}
\\
 I include below the theoretical results related to the fast extinction condition without demonstration in order to have a clear description of the SIS model \cite{Krivelevich1}

\begin{definition} \label{def1}
Define $S$ to be the $N \times N$ system matrix:\\
\begin{equation*}
S_{ij}= \left \{
\begin{array}{ll}
	1-\delta_{i}  &~~~~\text{if~} i=j \\
	r_{j} \beta_{ji} \frac{\gamma_{i}}{\gamma_{i}+\delta_{i}} & ~~~~\mbox{otherwise}
\end{array}
\right .
\end{equation*}
Let $|\lambda_{1,S}|$ be the magnitude of the largest eigenvalue and $\widehat{C}(t)=\sum_{i=1}^{N} p_{i}(t)$ the expected number of carriers at $t$ of the dynamical system.
\end{definition}

\begin{theorem} \label{TeoFastChak}
{\bf (Condition for fast extinction)}. Define $s=|\lambda_{1,S}|$ to be the {\bf survivability score} for the system. If $s=|\lambda_{1,S}|<1$, then we have fast extinction in the dynamical system, that is, $\widehat{C}(t)$ decays exponentially quickly over time.
\end{theorem}  
Where $|\lambda_{1,S}|$ is the magnitude of the largest eigenvalue of $S$, being $S$ an $N \times N$ system matrix defined as
$S_{ij}=1-\delta_{i}$ if $i=j$ and $S_{ij}=r_{j}\beta_{ji} \frac{\gamma_{i}}{\gamma_{i}+\delta_{i}}$ otherwise, and being $\widehat{C}(t)=\sum_{i=1}^{N}p_{i}(t)$ the expected number of carriers at time $t$ of the dynamical system.
Two additional results that appears in \cite{Deepa1} are the following

\begin{corollary}\label{CoroFastChak}
{\bf (Condition for fast extinction homogeneous case for Chakrabati SIS model)}
If $\delta_{i}= \delta, r_{i}=r, \gamma_{i}=\gamma$, for all $i$, and $B=[\beta_{ij}]$ is a symmetric binary matrix (links are undirected, and are always up or always down), then the condition for fast extinction is 
$\frac{\gamma}{\delta(\gamma+\delta)}\lambda_{1,B} < 1$
\end{corollary}

From the previous results \ref{TeoFastChak} and \ref{CoroFastChak} the fast extinction condition depends on the parameters $\delta, \gamma$ and $\beta$ as well as on the largest eigen-value 
$\lambda_{1,S}$ of the matrix of the dynamical system and therefore on the topology of the net.

Based on the equations \ref{probp} and \ref{probq} of the SIS model, implement a simulation in MATLAB whose results I show below. 

\begin{figure}[H]
\centering
\includegraphics[width=15cm, height=8cm]{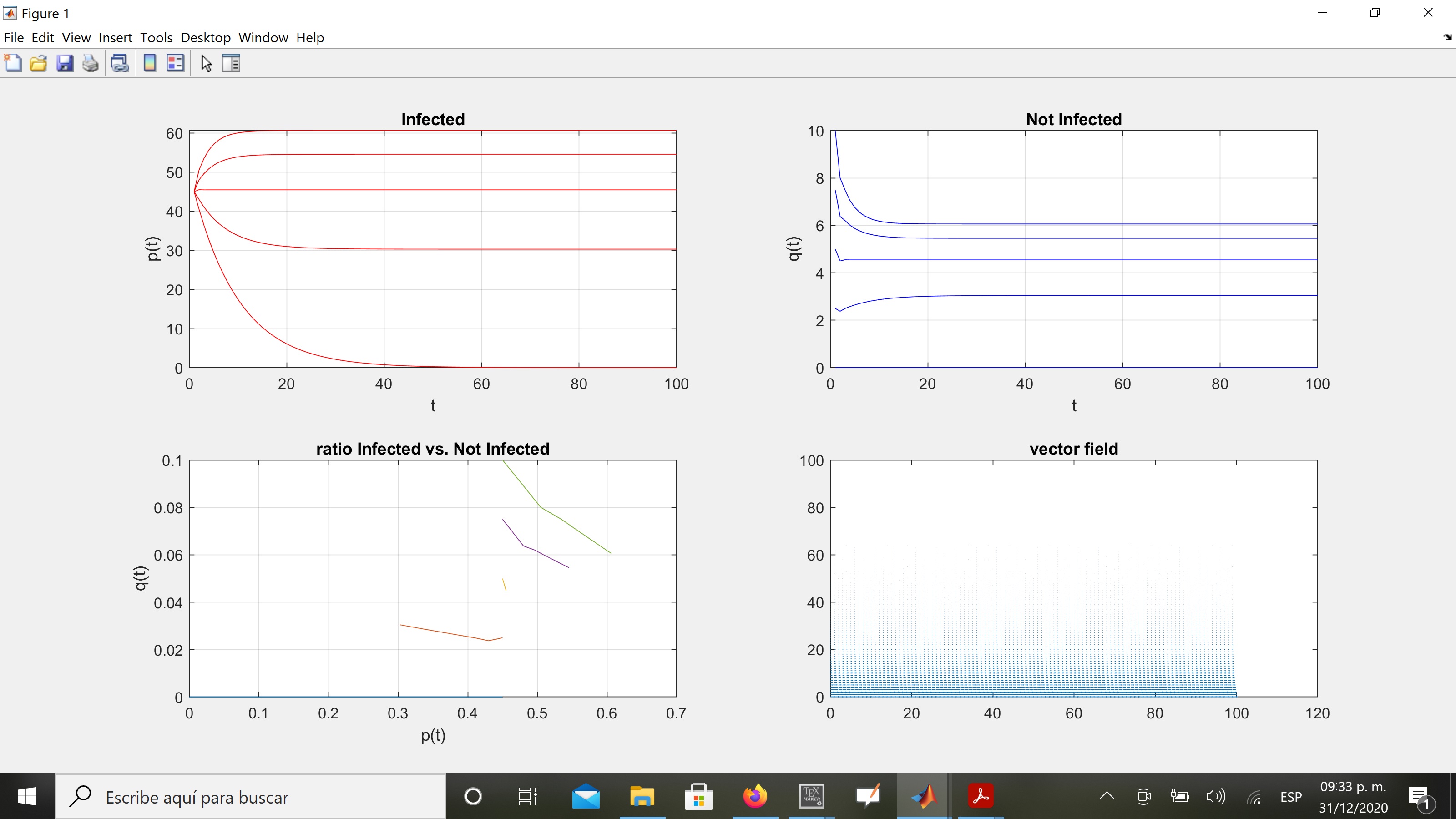}
%\centerline{\psfig{file=chakra101b.eps,width=4.2cm}}
%\centerline{\psfig{file=..\FigurasArticulo\MSEIR.jpeg,width=4.2cm}}
\vspace*{3pt}
\caption{Simulation of Chakrabarti SIS model with $\delta=0.1$, $\gamma=0.1+\Delta,\beta=0.1+\Delta$ on a Power law network}
\label{SISChak1}
\end{figure}

The figure \ref{SISChak1} shows the behavior of the Chakrabarti SIS model in a power law type network from five different initial conditions that consist in varying $\gamma$ and $\beta$ in increments of $0.05$.
Each process converges to states where the number of infected is higher than the number of susceptible.

\begin{figure}[H]
\centering
\includegraphics[width=15cm, height=8cm]{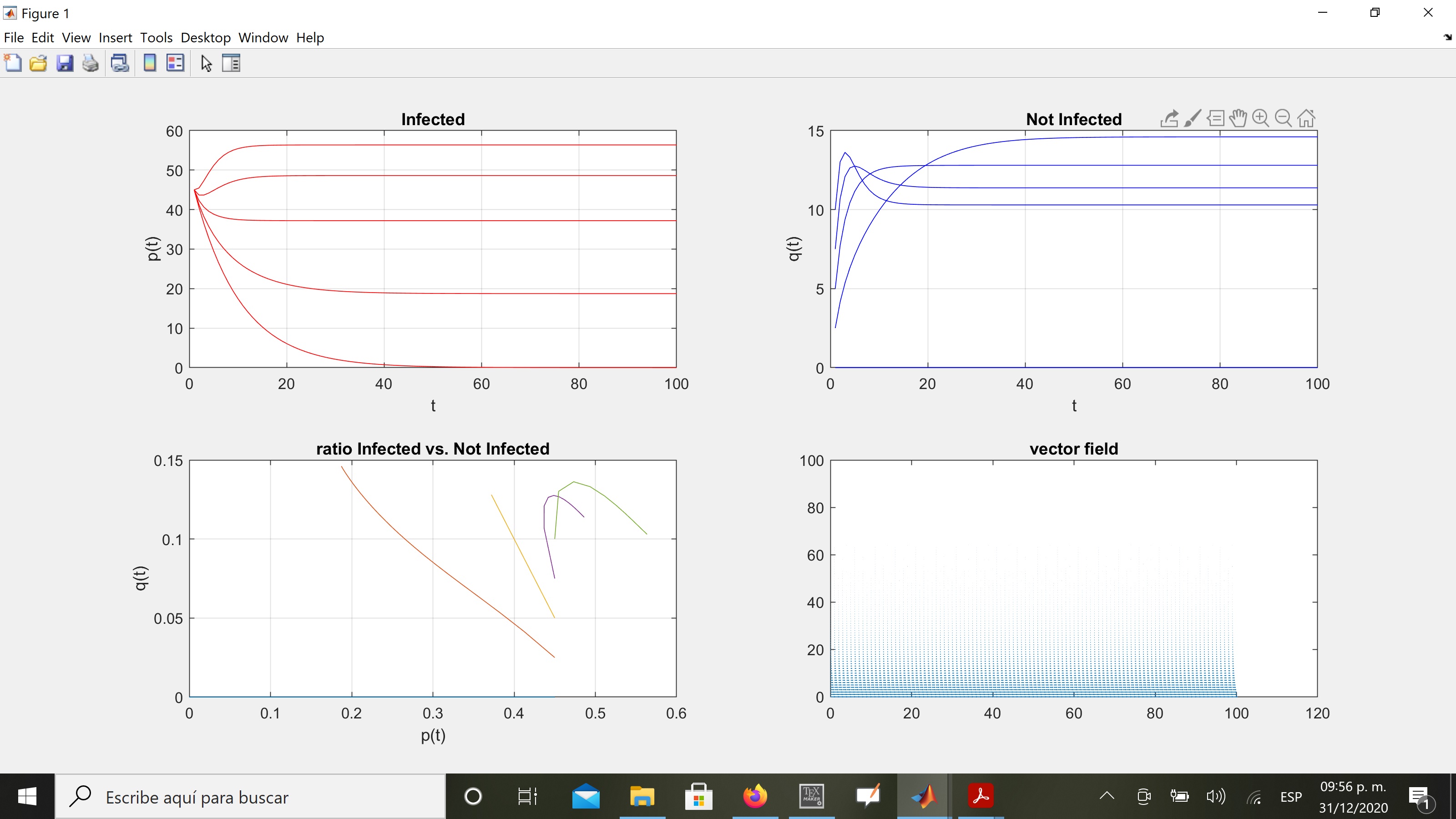}
%\centerline{\psfig{file=chakra101b.eps,width=4.2cm}}
%\centerline{\psfig{file=..\FigurasArticulo\MSEIR.jpeg,width=4.2cm}}
\vspace*{3pt}
\caption{Simulation of Chakrabarti SIS model with $\delta=0.1$, $\gamma=0.1+\Delta,\beta=0.1+\Delta$ on a Lattice 4 network}
\label{SISChak4}
\end{figure}

The figure \ref{SISChak4} shows the behavior of the Chakrabarti under the same conditions as the
behaviour of the simulation in figure \ref{SISChak1} but the topology of the network is a lattice.
It should be pointed out that process in figure \ref{SISChak4} have some variations with respect
to the one shown in figure \ref{SISChak1} due to the difference in the degree characteristics of the
respective topologies. 

\begin{figure}[H]
\centering
\includegraphics[width=15cm, height=8cm]{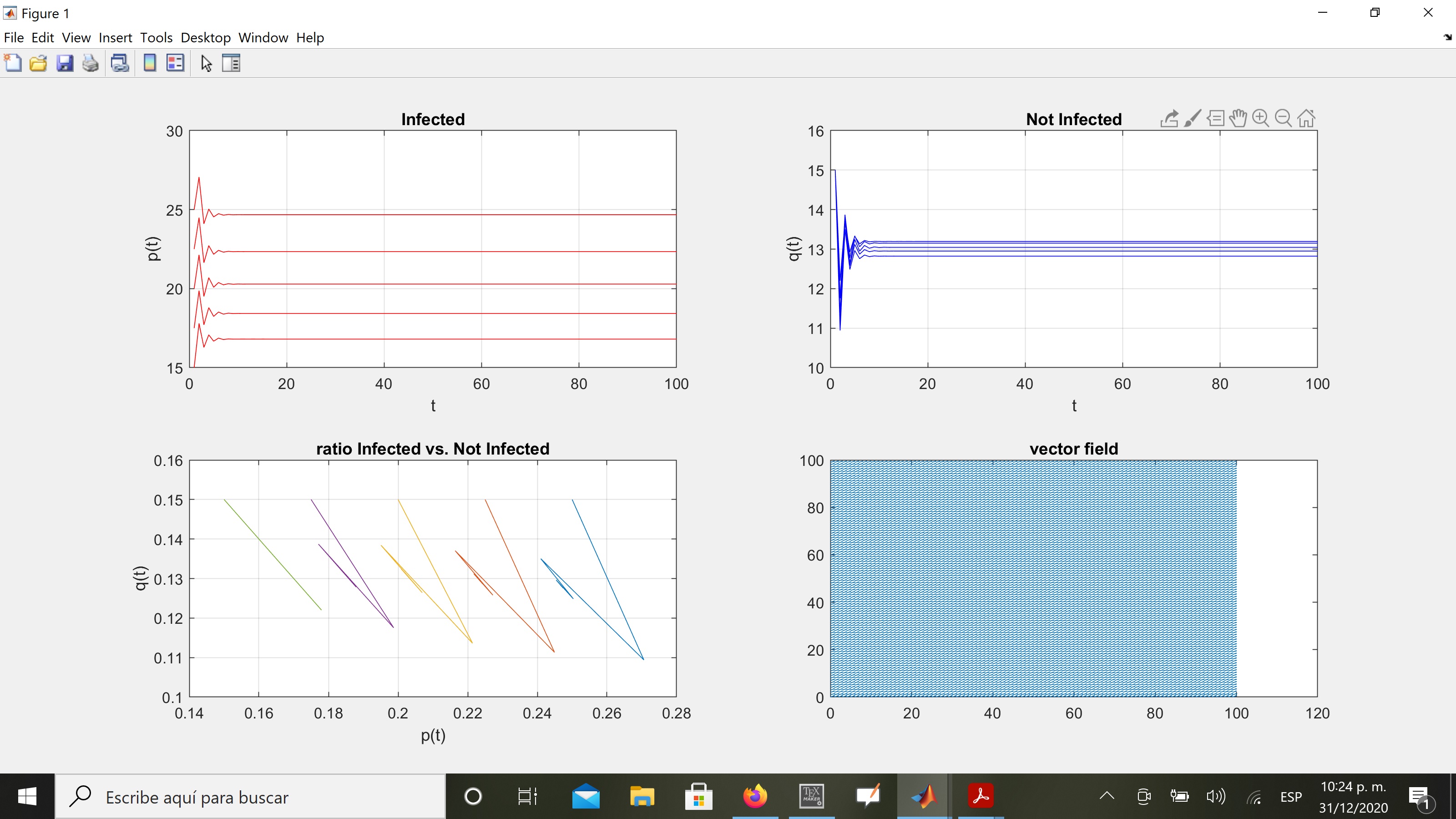}
%\centerline{\psfig{file=chakra101b.eps,width=4.2cm}}
%\centerline{\psfig{file=..\FigurasArticulo\MSEIR.jpeg,width=4.2cm}}
\vspace*{3pt}
\caption{Simulation of Chakrabarti SIS model with $\delta=0.5+\Delta$, $\gamma=0.3,\beta=0.4$ on a Power law network}
\label{SISChak1b}
\end{figure}
The figure \ref{SISChak1b} shows the behavior of the Chakrabarti SIS model in a power law type network from five different initial conditions that consist in varying $\gamma$ in increments of $0.05$.
Each process converges to states where the number of infected is higher than the number of susceptible.

\begin{figure}[H]
\centering
\includegraphics[width=15cm, height=8cm]{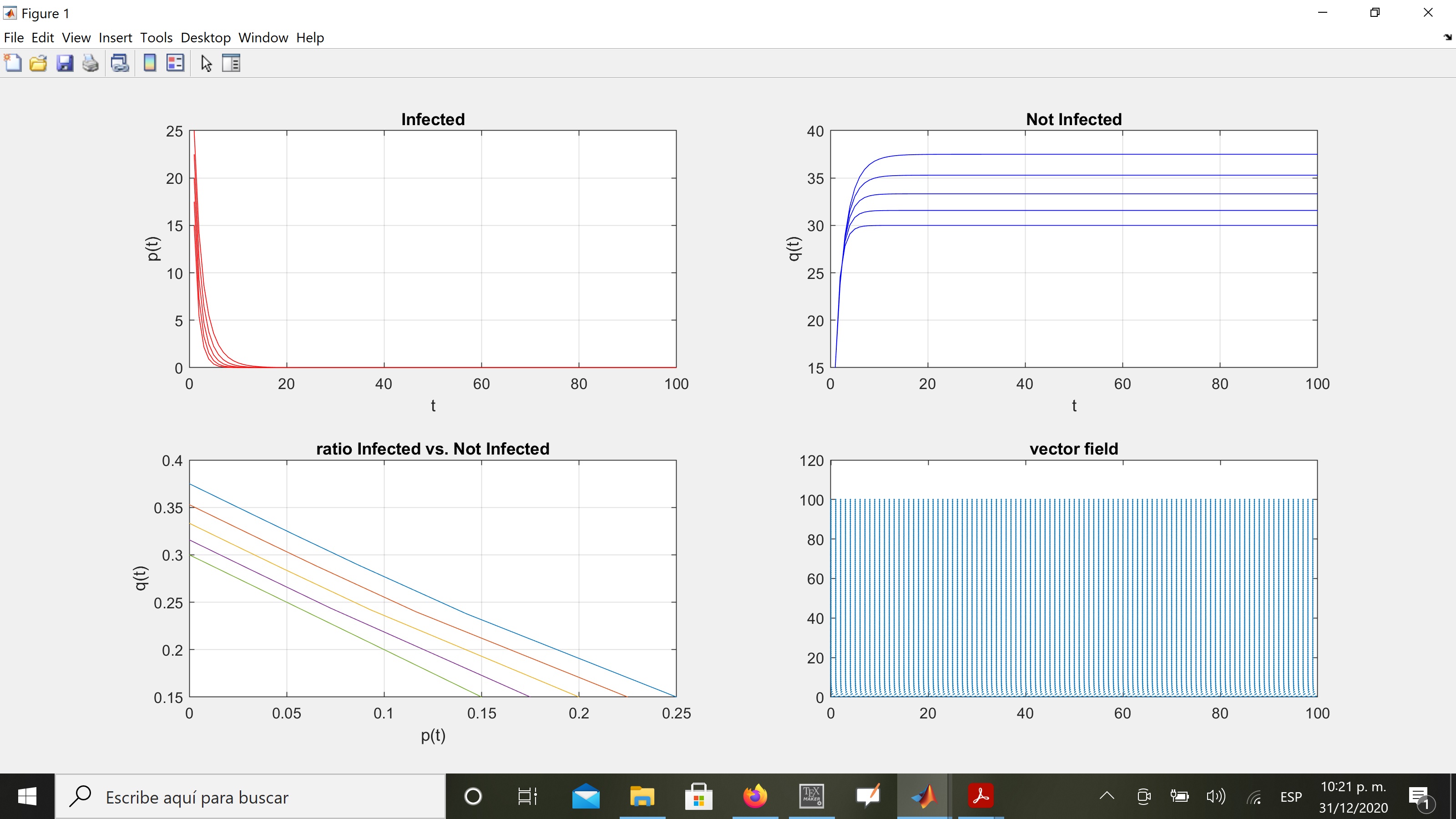}
%\centerline{\psfig{file=chakra101b.eps,width=4.2cm}}
%\centerline{\psfig{file=..\FigurasArticulo\MSEIR.jpeg,width=4.2cm}}
\vspace*{3pt}
\caption{Simulation of Chakrabarti SIS model with $\delta=0.5+\Delta$, $\gamma=0.3,\beta=0.4$ on a Lattice 4 network}
\label{SISChak4b}
\end{figure}
The figure \ref{SISChak4b} shows the behavior of the Chakrabarti SIS model in a Lattice 4 type network from five different initial conditions that consist in varying $\gamma$ in increments of $0.05$.
It can be seen that the process achieves fast extinction condition given that the topology ia lattice. So with this example we can observe that the topology can make a difference between having a fast extinction of the epidemic or converging to a state where the number of infected is bigger the number of suceptible and stays stable in such state.

In \cite{CRL1} it was formulated a similar discrete epidemic model proposed by the authors of \cite{Deepa1} having one additional states in order to let the nodes to be warned by a message or 
to receive a vaccine. 

The idea behind this additional states was to explore prevention alternatives as well as the possible eradication of a virus in a computer network either through warning messages or by distribution of a vaccine \cite{chayes1}.
The results regarding the fast extinction condition of the virus as well as the fixed point results were very similar to those of the \cite{Deepa1} model.

\begin{figure}[H] 
\centering
\includegraphics[width=10cm, height=6cm]{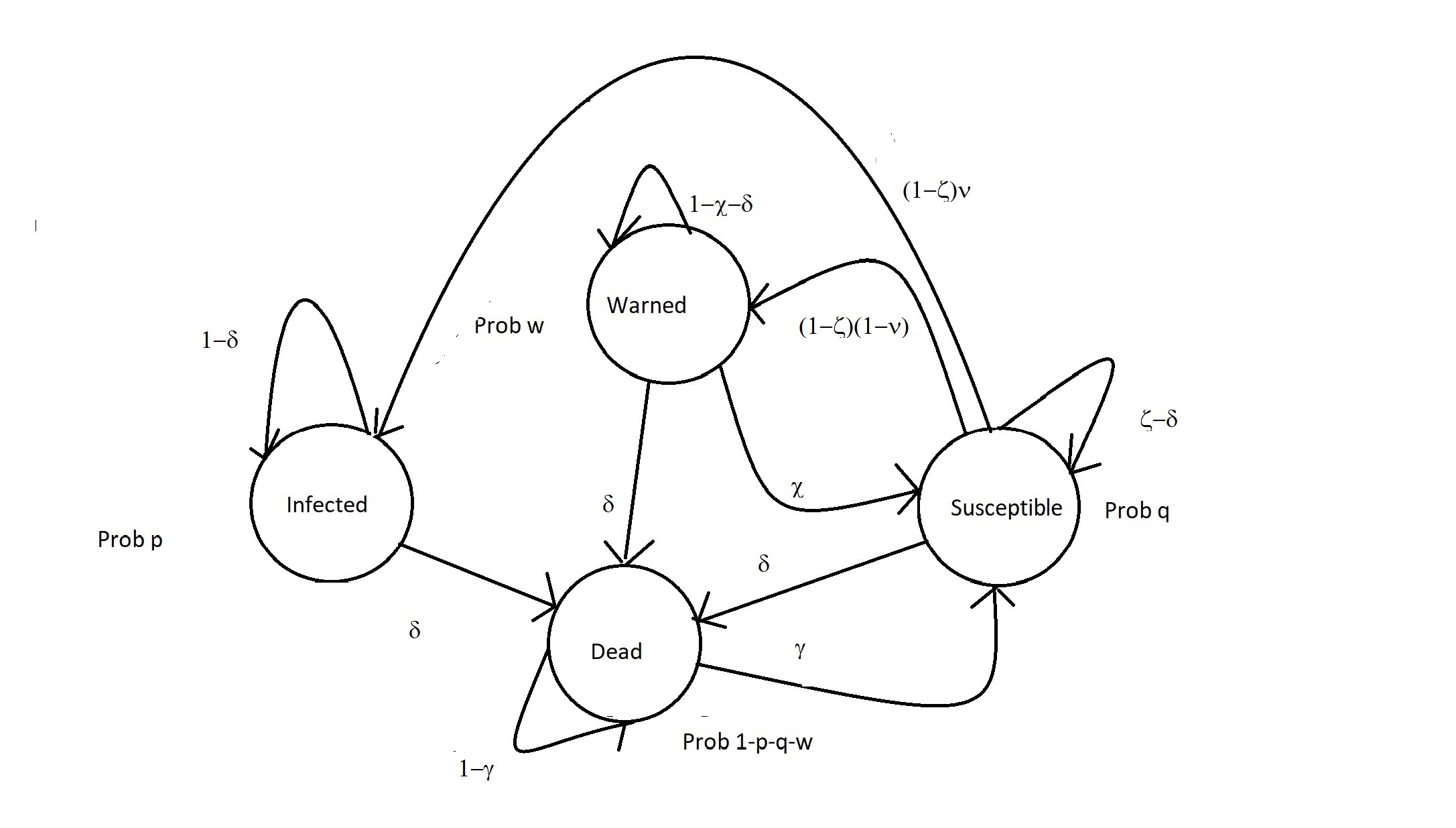}
%\centerline{\psfig{file=chakra101b.eps,width=4.2cm}}
%\centerline{\psfig{file=..\FigurasArticulo\MSEIR.jpeg,width=4.2cm}}
\vspace*{5pt}
\caption{The SIRS model}
\label{MiSIRS1fig}
\end{figure}

Each node can be in one of three states: \emph{Infected},\emph{Warn Info}, \emph{No Info} or \emph{Dead}, with transitions between them as shown in Diagram \ref{MiSIRS1fig}.
%Our model can be graphically represented as follows:\\
The next graph represent the transitions that take place in each node for this model.\\

%\begin{figure}[th]

%\subsection{Fast extinction result} \label{fastmimodelo}

Making the same node independence probability assumption that is stated in equation (\ref{indep}) and
taking into account the new states and transition probabilities shown in the figure \ref{MiSIRS1fig}, the equations (\ref{probp}) and (\ref{probq}) as well as the new equation corresponding to $w_{i}$ 
%and $h_{i}$
can be expressed as follows:
\begin{eqnarray}
p_{i}(t)&=&p_{i}(t-1)(1-\delta_{i})+q_{i}(t-1)(1-\zeta_{i}(t))\nu_{i} \label{probnp}\\
 &  & \nonumber\\
q_{i}(t)&=&q_{i}(t-1)(\zeta_{i}(t)-\delta_{i})+ \label{probnq}\\
	&  &(1-p_{i}(t-1)-q_{i}(t-1)-w_{i}(t-1))\gamma_{i} \nonumber\\
	&  &+\chi_{i}w_{i}(t-1) \nonumber \\
 &  & \nonumber\\
w_{i}(t)&=&(1-\zeta_{i}(t))(1-\nu_{i})q_{i}(t-1) \label{probw}  \\
    &  & +(1-\chi_{i}-\delta_{i})w_{i}(t-1) \nonumber 
\end{eqnarray}
%%%%%%
Based on the equations \ref{probnp},\ref{probq} and  \ref{probw} of the SIRS model, I implemented a simulation in MATLAB whose results I show below. 

\begin{figure}[H]
\centering
\includegraphics[width=15cm, height=8cm]{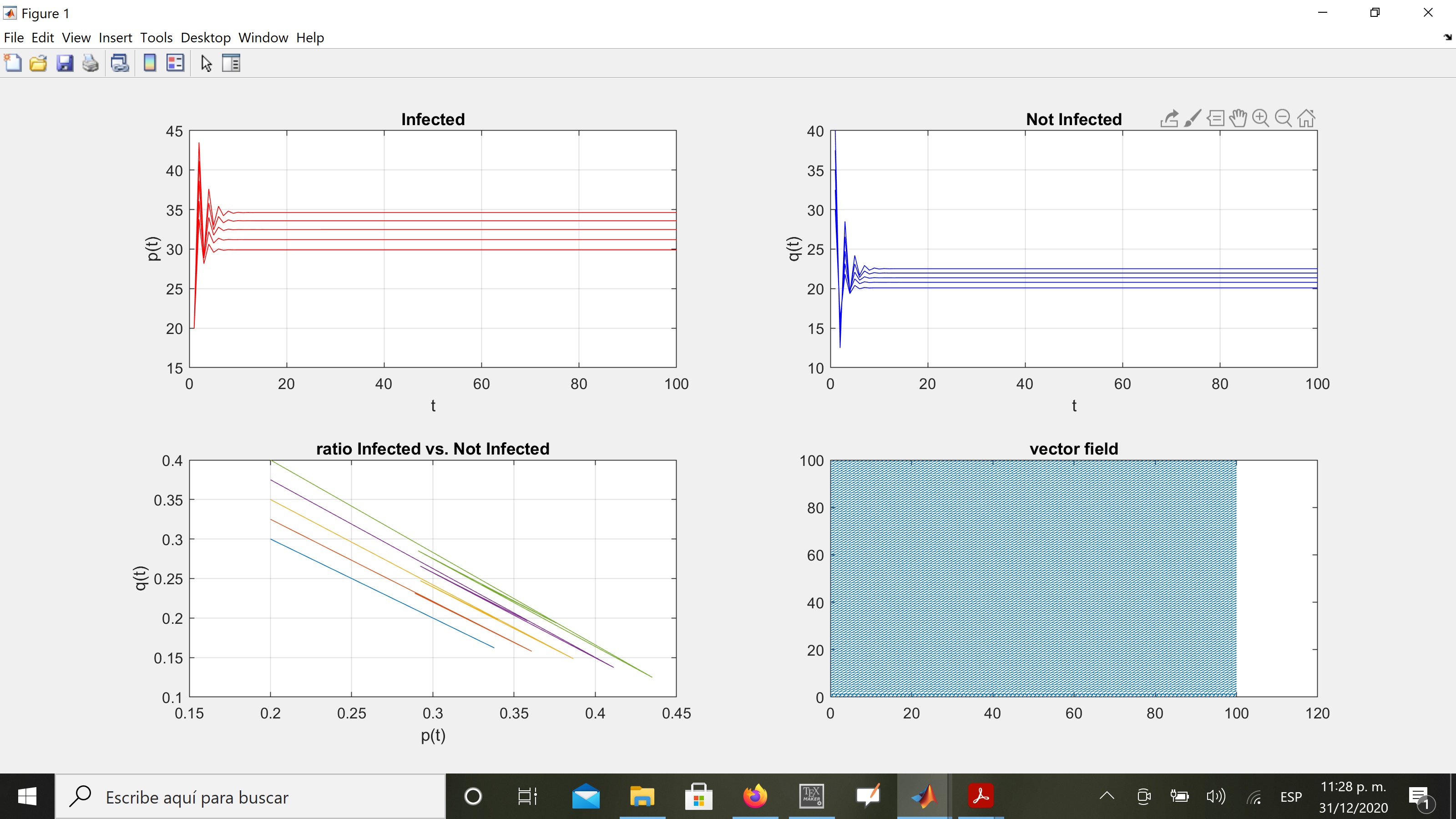}
%\centerline{\psfig{file=chakra101b.eps,width=4.2cm}}
%\centerline{\psfig{file=..\FigurasArticulo\MSEIR.jpeg,width=4.2cm}}
\vspace*{3pt}
\caption{Simulation of SIRS model with $\delta=0.6$, $\gamma=0.6+\Delta, \nu=1,\chi=1,\beta=0.3$ on a Power law network}
\label{SIRS1}
\end{figure}

The figure \ref{SIRS1} shows the behavior of the SIRS model in a power law type network from five different initial conditions that consist in varying $\gamma$ and $\beta$ in increments of $0.05$.
Each process converges to states where the number of infected is higher than the number of susceptible.

\begin{figure}[H]
\centering
\includegraphics[width=15cm, height=8cm]{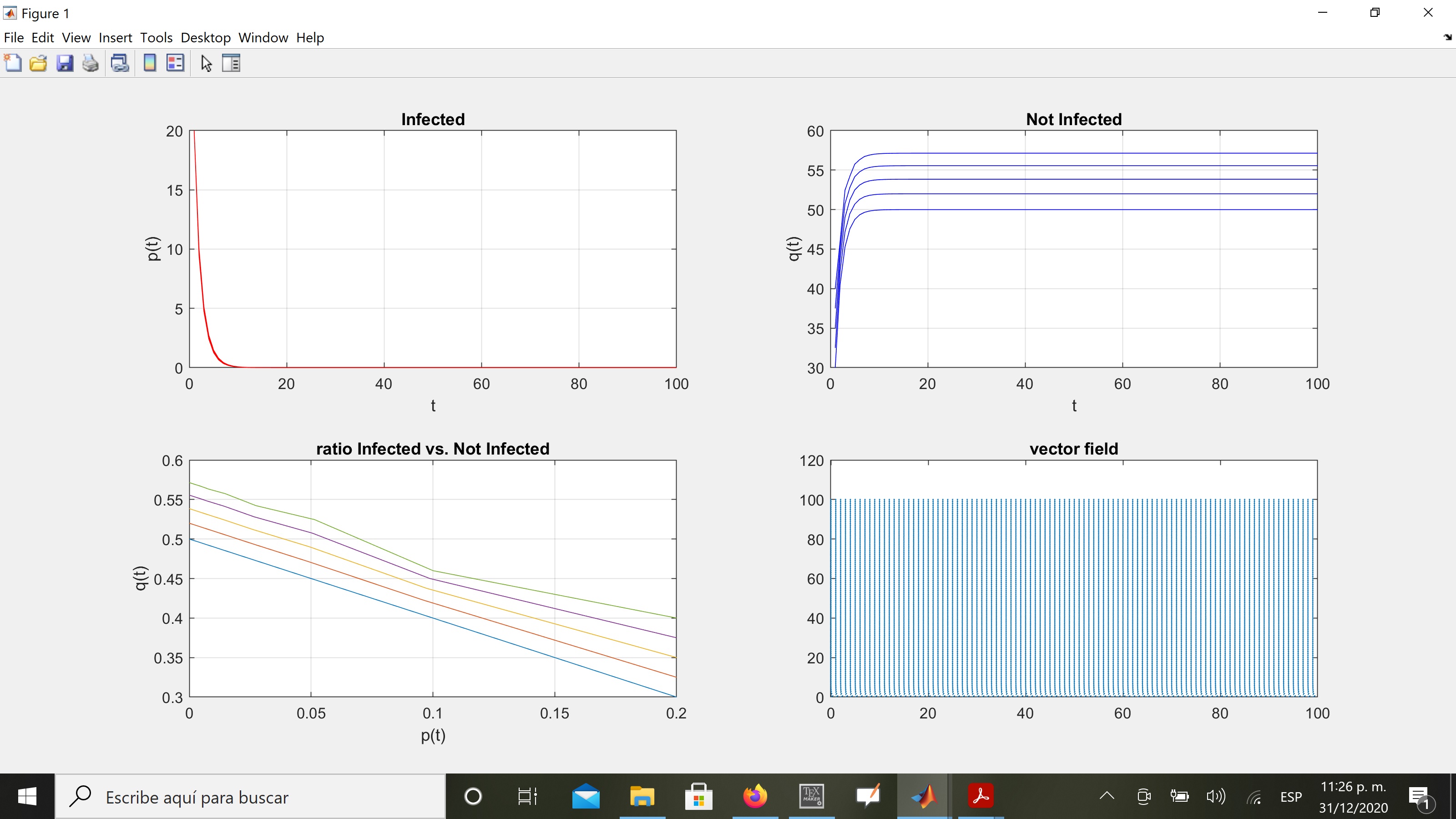}
%\centerline{\psfig{file=chakra101b.eps,width=4.2cm}}
%\centerline{\psfig{file=..\FigurasArticulo\MSEIR.jpeg,width=4.2cm}}
\vspace*{3pt}
\caption{Simulation of SIRS model with $\delta=0.6$, $\gamma=0.6+\Delta, \nu=1,\chi=1,\beta=0.3$ on a Lattice 4 network}
\label{SIRS4}
\end{figure}

From the behaviour shown in figure \ref{SIRS1} and \ref{SIRS4} of the SIRS we can claim that
the topology has an impact for achieving the fast extinction of an epidemic.

%% Hasta aqui voy CRL 31/dic/2020

\section{Ideas about isolation strategies} \label{isolation}
As can be seen in the results obtained in the simulations of both the Chakrabarti SIS model and the SIRS model, the network topology has a significant impact on the rapid extinction of a computer virus in a network or its permanence in it. Given that the topology of the network is related to its adjacency matrix, a possible isolation strategy to minimize the first eigenvalue of said matrix \cite{Van1}. The question is that this strategy would be reduced to the search for a Hamiltonian circuit, which in general is an NP-hard problem. 
One possibility is to look for algorithms that approximate the obtaining of the minimum value Hamiltonian path by algorithmic techniques of the closest neighbors type that seek to approximate the minimum value Hamiltonian cycle with some guarantee of approximation
and eliminate from the original interconnection graph those edges 
that do not belong to said cycle. This is what is proposed in \cite{Van2}. 

Given that the problem of minimization of eigenvalues of the adjacency matrix is NP-hard and the algorithms of approximation to the problem of calculation of Hamiltonian cycles that 
guarantee a certain quality of approximation, I would propose to try to transform 
any power law graph to a graph of type lattice 4 to reduce the mentioned eigenvalue 
in some way, even if it is not the minimum.
An application to contain the spread of a virus in a network, be it computerized or of people, could be to detect the value of the parameters that characterize a process of diffusion of a virus taking into account the interconnection medium where the epidemic process will take place and from there determine which network nodes to isolate by modifying the adyascences of the associated graph in such a way as to obtain a notable reduction in the spread of a virus or even achieve its rapid extinction.

\section{Conclusions and Future Work}
In this article, he makes a historical account of the classic mathematical models for the study of epidemic processes that made it possible to develop models to address problems of virus spread in networks. We were also able to observe that the concerns that arise when trying to solve problems that arise in the field of virus spread can contribute ideas to the development of mathematical models in epidemiology by incorporating aspects of the interconnection networks on which epidemics take place. By introducing the structures of interconnection networks to epidemiological models, we can obtain elements that guide epidemiologists in making decisions about isolation strategies to try to contain and eventually eradicate a disease through the development of strategies based on knowledge. of the dissemination of information in interconnection networks.

% hasta aqui voy CRL 22/mayo/2020

%\bibliographystyle{ws-acs}

%\bibliography{carlos}
%\bibliographystyle{plain}

\end{document}